\documentclass[a4paper,12pt]{amsart}
\usepackage{amsfonts}
\usepackage{amssymb}
\usepackage{ifthen}
\usepackage{graphicx}
\usepackage{enumitem}
\usepackage{amsmath}
\nonstopmode \numberwithin{equation}{section}
\usepackage[usenames]{color}
\newtheorem{thm}{Theorem}[section]
\newtheorem{lem}{Lemma}[section]
\newtheorem{cor}[thm]{Corollary}
\newtheorem{prop}[thm]{Proposition}

\theoremstyle{definition}
\newtheorem{mlem}{Main lemma}[section]
\newtheorem{assertion}{Assertion}[section]
\newtheorem{cl}{Claim}[section]
\newtheorem{ca}{Case}[section]
\newtheorem{sca}{Subcase}[section]
\newtheorem{scl}{Subclaim}[section]
\newtheorem{conj}[thm]{Conjecture}
\newtheorem{fact}{Fact}[section]
\newtheorem{defn}{Definition}[section]
\newtheorem{prob}{Problem}[section]
\newtheorem{ques}[thm]{Question}
\newtheorem{rem}{Remark}[section]
\newtheorem{exam}{Example}[section]

\numberwithin{equation}{section}

\newcounter {own}
\def\theown {\thesection       .\arabic{own}}

\newlist{steps}{enumerate}{1}
\setlist[steps,1]{
  leftmargin=*,
  label=\textbf{Step \arabic*}.,
  ref=Step~\arabic*,
}

\newenvironment{pf}[1][]{%
 \vskip 3mm
 \noindent
 \ifthenelse{\equal{#1}{}}%
  {{\slshape Proof. }}%
  {{\slshape #1.} }%
 }%
{\qed\bigskip}

\newcounter{alphabet}
\newcounter{tmp}
\newenvironment{Thm}[1][]{\refstepcounter{alphabet}%
\bigskip%
\noindent%
{\bf Theorem \Alph{alphabet}}%
\ifthenelse{\equal{#1}{}}{}{ (#1)}%
{\bf .} \itshape}{\vskip 8pt}

\newcounter{alphabet2}

\makeatletter
\newcommand{\Ref}[1]{\@ifundefined{r@#1}{}{\setcounter{tmp}{\ref{#1}}\Alph{tmp}}}
\makeatother

\newenvironment{Lem}[1][]{\refstepcounter{alphabet}%
\bigskip%
\noindent%
{\bf Lemma \Alph{alphabet}}%
{\bf .} \itshape}{\vskip 8pt}

\newcommand{\IC}{{\mathbb C}}
\newcommand{\ID}{{\mathbb D}}

\def\be{\begin{equation}}
\def\ee{\end{equation}}

\newcommand{\ben}{\begin{enumerate}}
\newcommand{\een}{\end{enumerate}}

\newcommand{\blem}{\begin{lem}}
\newcommand{\elem}{\end{lem}}
\newcommand{\bthm}{\begin{thm}}
\newcommand{\ethm}{\end{thm}}
\newcommand{\bcor}{\begin{cor}}
\newcommand{\ecor}{\end{cor}}
\newcommand{\beg}{\begin{exam}}
\newcommand{\eeg}{\end{exam}}
\newcommand{\begs}{\begin{examples}}
\newcommand{\eegs}{\end{examples}}
\newcommand{\bdefe}{\begin{defn}}
\newcommand{\edefe}{\end{defn}}
\newcommand{\bprob}{\begin{prob}}
\newcommand{\eprob}{\end{prob}}
\newcommand{\bques}{\begin{ques}}
\newcommand{\eques}{\end{ques}}
\newcommand{\bei}{\begin{itemize}}
\newcommand{\eei}{\end{itemize}}
\newcommand{\bcon}{\begin{conj}}
\newcommand{\econ}{\end{conj}}
\newcommand{\bop}{\begin{op}}
\newcommand{\eop}{\end{op}}

\newcommand{\bas}{\begin{assertion}}
\newcommand{\eas}{\end{assertion}}

\newcommand{\bfa}{\begin{fact}}
\newcommand{\efa}{\end{fact}}

\newcommand{\bca}{\begin{ca}}
\newcommand{\eca}{\end{ca}}
\newcommand{\bsca}{\begin{sca}}
\newcommand{\esca}{\end{sca}}

\newcommand{\bcl}{\begin{cl}}
\newcommand{\ecl}{\end{cl}}

\newcommand{\bmlem}{\begin{mlem}}
\newcommand{\emlem}{\end{mlem}}

\newcommand{\bscl}{\begin{scl}}
\newcommand{\escl}{\end{scl}}

\newcommand{\bcons}{\begin{conjs}}
\newcommand{\econs}{\end{conjs}}

\newcommand{\bprop}{\begin{prop}}
\newcommand{\eprop}{\end{prop}}

\newcommand{\br}{\begin{rem}}
\newcommand{\er}{\end{rem}}
\newcommand{\brs}{\begin{rems}}
\newcommand{\ers}{\end{rems}}
\newcommand{\bo}{\begin{obser}}
\newcommand{\eo}{\end{obser}}
\newcommand{\bos}{\begin{obsers}}
\newcommand{\eos}{\end{obsers}}
\newcommand{\bpf}{\begin{pf}}
\newcommand{\epf}{\end{pf}}
\newcommand{\ba}{\begin{array}}
\newcommand{\ea}{\end{array}}
\newcommand{\beq}{\begin{eqnarray}}
\newcommand{\beqq}{\begin{eqnarray*}}
\newcommand{\eeq}{\end{eqnarray}}
\newcommand{\eeqq}{\end{eqnarray*}}

\newcounter{minutes}
\divide\time by 60
\newcounter{hours}
\multiply\time by 60 \addtocounter{minutes}{-\time}

\begin{document}
\bibliographystyle{amsplain}
\title{Boundary Schwarz lemma for harmonic mappings having zero of order $p$}
\thanks{
\indent
File:~\jobname .tex,
printed: \number\day-\number\month-\number\year,
\thehours.\ifnum\theminutes<10{0}\fi\theminutes
}
\author{Xiao-Jin Bai}
\address{Xiao-Jin Bai, School of Mathematical Sciences, Huaqiao
University, Quanzhou 362021, China.}
\email{xiaojin\_bai@foxmail.com}

\author{Jie Huang}
\address{Jie Huang, School of Mathematical Sciences, Huaqiao
	University, Quanzhou 362021, China and Technion-Israel institute of Technology, Guangdong Technion, 241 Daxue Road, Shantou,  Guangdong 515063, China.}
\email{jie.huang@gtiit.edu.cn;}

\author{Jian-Feng Zhu}
\address{Jian-Feng Zhu, School of Mathematical Sciences, Huaqiao
University, Quanzhou 362021, China.}
\email{flandy@hqu.edu.cn}


\begin{abstract}
Suppose $w$ is a sense-preserving harmonic mapping of the unit disk $\ID$ such that $w(\ID)\subseteq\ID$ and
$w$ has a zero of order $p\geq1$ at $z=0$. In this paper, we first improve
the Schwarz lemma for $w$, and then, we establish its boundary Schwarz lemma.
Moreover, by using the automorphism of $\ID$, we further generalize this result.

\end{abstract}

\thanks{
}

\maketitle \pagestyle{myheadings}
\markboth{Bai, Huang and Zhu }{Boundary Schwarz lemma for harmonic mappings having zero of order $p$}

\section{Introduction}
Let $\ID=\{z\in\mathbb{C}:\; |z|<1\}$ be the unit disk, $\mathbb{T}=\{z\in\mathbb{C}:\; |z|=1\}$ the unit circle, and $\overline{\ID}$ the closure of $\mathbb{D}$, i.e., $\overline{\ID}=\ID\cup \mathbb{T}$.
For $z\in \ID$, the {\it formal derivatives} of a complex-valued function $f$ are defined by:
\beqq\label{eq1.1}
f_z=\frac{1}{2}\left(\frac{\partial f}{\partial x}-i\frac{\partial f}{\partial y}\right)\;\;\mbox{and}\;\;f_{\bar{z}}=\frac{1}{2}\left(\frac{\partial f}{\partial x}+i\frac{\partial f}{\partial y}\right).
\eeqq
For each $\alpha\in[0,2\pi]$, the {\it directional derivative} of $f$ at $z\in \ID$ is defined by
\beqq\label{eq1.2}
\partial_{\alpha}f(z)=\lim\limits_{r\rightarrow0^+}\frac{f(z+re^{i\alpha})-f(z)}{r}=e^{i\alpha}f_z(z)+e^{-i\alpha}f_{\bar{z}}(z).
\eeqq
Then
\beqq\label{eq1.3}
\max\limits_{0\leq\alpha\leq2\pi}\{|\partial_{\alpha}f(z)|\}=\Lambda_f(z)=|f_z(z)|+|f_{\bar{z}}(z)|
\eeqq
and
\beqq\label{eq1.4}
\min\limits_{0\leq\alpha\leq2\pi}\{|\partial_{\alpha}f(z)|\}=\lambda_f(z)=\big||f_z(z)|-|f_{\bar{z}}(z)|\big|.
\eeqq
A function $f$ is said to be locally univalent and sense-preserving in $\ID$ if and only if its Jacobian $J_f$ satisfies the following condition
(cf. \cite{HLew1936}): For any $z \in \ID$,
$$J_f(z)=|f_z(z)|^2-|f_{\bar{z}}(z)|^2>0.$$

Here and hereafter, the notation $C^m(E)$ denotes the set of all functions which are $m$-times continuously differentiable in domain $E\subset \IC$, where $m\geq 0$ is an integer. In particular, $C^0(E)$, which is always denoted by $C(E)$, means the set of all continuous functions in $E$.

A function $w\in C^2(E)$ is said to be harmonic in $E$ if it satisfies the following Laplace equation
$$\Delta w=4w_{z\bar{z}}=0.$$
Obviously, harmonic mappings are generalizations of analytic functions.

In a simply connected domain $\Omega\subset\IC$, a harmonic mapping $w$ has the representation $w=h+\bar{g}$, where
$h$ and $g$ are analytic in $\Omega$. Furthermore, if $g(0)=0$, then the representation is unique and called the {\it canonical representation}.
We refer to
\cite{Du-04} for more properties of harmonic mappings.

In the rest of this paper, we use $w$ to stand for the harmonic mappings of $\ID$, and $f$ to stand for the analytic function of $\ID$.

\subsection{The multiplicity of zeros for analytic functions and harmonic mappings}
\subsubsection{Analytic case.}
Suppose that $f$ is an analytic function of $\ID$. Then $f$ is said to have a zero of order $n$ at $z_0$, where $n\geq1$, denoted by $\mu(z_0, f)=n$,
if $f(z_0)=Df(z_0)=\cdots=D^{n-1}f(z_0)=0$ and $D^nf(z_0)\neq0$, i.e.,
$$f(z)=\sum\limits_{k=n}^{\infty}a_k(z-z_0)^k,\ \ \ \mbox{for}\  \  z\in\ID.$$
Here and hereafter the symbol $D^kf$ (resp. $\bar{D}^k f$) means the $k-$th order derivative with respect to $z$ (resp. $\bar{z}$) of the complex-valued function $f$, i.e., $D^kf=(\frac{\partial}{\partial z})^kf(z)$ (resp. $\bar{D}^kf=(\frac{\partial}{\partial \bar{z}})^kf(z)$).

The following result is a consequence of the Schwarz-Pick lemma applied to the function $f/z^p$ (cf. \cite[Corollary 1.3]{JBG1981} or \cite[Remark 3]{Osserman00}).

\begin{Lem}\label{HLZ-lemma4}
Let $f:\mathbb{D} \to \mathbb{D}$ be an analytic function with $\mu (0,f) = p \ge 1$. Then
for any $z\in\ID$,
$$|f(z)|\leq |z|^p\frac{|z|+|a_p|}{1+|a_p||z|},$$
where $a_p=\frac{D^pf(0)}{p!}$.
\end{Lem}

\subsubsection{Harmonic case.}
Suppose that $w=h+\bar{g}$ is a harmonic mapping of $\ID$.
For any $z\in\ID$, let
$$\omega(z)=\frac{\overline{w_{\bar{z}}(z)}}{w_z(z)}$$
be the {\it second complex dilatation} of $w$. Then
$\omega(z)=\frac{g'(z)}{h'(z)}$
is an analytic function in $\ID$.
Moreover, if $w(z)$ is sense-preserving, then $|\omega(z)|<1$ for all $z\in\ID$.

We now introduce the definition of the multiplicity for sense-preserving harmonic mappings $w$ in $\ID$.
Suppose that $w=h+\bar{g}$ is a sense-preserving harmonic mapping of $\ID$, where
$h$ and $g$ have respectively multiplicity $n$ and $m$ at $z_0$ with $w(z_0)=0$, i.e.,
$$h(z)=\sum\limits_{k=n}^\infty a_k(z-z_0)^k,\ \ g(z)=\sum\limits_{k=m}^\infty b_k(z-z_0)^k,\ \ \ z\in\ID.$$
Then $n<m$ or $m=n$ and $|b_n|<|a_n|$, since $|\omega(z_0)|<1$.
We say that $w$ has a zero of order $n$ at $z_0$ and write $\mu(z_0, w)=n$.

The following lemma is due to Ponnusamy and Rasila \cite{Samy-Antti}.
Note that if $p=1$, then it is the well-known harmonic version of the classical Schwarz lemma due to Heinz \cite{Heinz}.

\begin{Lem}\label{Samy-Antti-lem}
Let $w$ be a sense-preserving harmonic mapping of $\ID$ such that $\mu(0, w)=p\geq1$ and $w(\ID)\subset\ID$. Then
for any $z\in\ID$,
$$
|w(z)|\leq\frac{4}{\pi}\arctan|z|^p\leq\frac{4}{\pi}|z|^p.
$$
\end{Lem}

Using Lemma \Ref{HLZ-lemma4}, we first improve Lemma \Ref{Samy-Antti-lem} as follows:

\blem\label{bhz-lem-1}
Let $w=h+\bar{g}$ be a sense-preserving harmonic mapping of $\ID$ such that $\mu(0, w)=p\geq1$ and $w(\ID)\subset\ID$. Then
for any $z\in\ID$,
\be\label{bhz-lem1-1}|w(z)|\le \frac{4}{\pi}\arctan \left[|z|^p\frac{|z|+\frac{\pi}{4}(|a_p|+|b_p|)}{1+\frac{\pi}{4}(|a_p|+|b_p|)|z|}\right],\ee
where $a_p=\frac{D^ph(0)}{p!}$ and $b_p=\frac{D^pg(0)}{p!}$.
\elem

Since $w$ is a harmonic self-mapping of $\ID$, it follows from \cite[Lemma 1]{Shaolin-Bloch} that
\be\label{April-22}|a_n|+|b_n|\leq\frac{4}{\pi},\ \ \  \mbox{for all}\ \  n=1, 2, \cdots. \ee
For any $0\leq r<1$, the function $\varphi(x)=\frac{r+\frac{\pi}{4}x}{1+\frac{\pi}{4}xr}$ is an increasing function of $x$, then we see that
$$\frac{4}{\pi}\arctan \left[|z|^p\frac{|z|+\frac{\pi}{4}(|a_p|+|b_p|)}{1+\frac{\pi}{4}(|a_p|+|b_p|)|z|}\right]\leq\frac{4}{\pi}\arctan |z|^p.$$

\subsection{The boundary Schwarz lemma for analytic functions and harmonic mappings}

Let us recall the following classical boundary Schwarz lemma for analytic functions, which was proved in \cite{JBG1981}.

\begin{Thm}\label{Thm-1}$($\cite[Page 42]{JBG1981}$)$
Suppose $f:\mathbb{D}\rightarrow\mathbb{D}$ is an analytic function with $f(0)=0$, and, further, $f$ is analytic
at $z=1$ with $f(1)=1$. Then, the following two conclusions hold:\ben
\item
 $f'(1)\geq1$.
 \item
 $f'(1)=1$ if and only if $f(z)\equiv z$.
 \een
\end{Thm}

Theorem \Ref{Thm-1} has the following generalization.

\begin{Thm}\label{Thm-2}$($\cite[Theorem 1.1$^{\prime}$]{LTPAMQ2015}$)$
Suppose $f:\mathbb{D}\rightarrow\mathbb{D}$ is an analytic function with $f(0)=0$, and, further, $f$ is analytic
at $z=\alpha\in\mathbb{T}$ with $f(\alpha)=\beta\in\mathbb{T}$. Then, the following two conclusions hold:
\ben
\item
 $\overline{\beta}f'(\alpha)\alpha\geq1$.
\item
$\overline{\beta}f'(\alpha)\alpha=1$ if and only if $f(z)\equiv e^{i\theta}z$, where $e^{i\theta}=\beta\alpha^{-1}$ and
$\theta\in\mathbb{R}$.
\een
\end{Thm}

We remark that, when $\alpha=\beta=1$, Theorem \Ref{Thm-2} coincides with Theorem \Ref{Thm-1}.

This useful result has attracted much attention and has been
generalized in various forms (see, e.g., \cite{BK94, CK, CP-0,
Kra11, LT16, Zhu}). Recently, Wang et. al. obtained the boundary Schwarz lemma for
solutions to the Poisson's equation (\cite{WZ18}).
By analogy with the studies in the above results, in this paper, we discuss the boundary Schwarz lemma for harmonic mappings having a zero of order $p$.
Our main results are as follows:

\begin{thm}\label{thm-1.1}
Let $w=h+\bar{g}$ be a sense-preserving harmonic mapping of $\ID$ such that $\mu(0, w)=p\geq1$ and $w(\ID)\subset\ID$.
If $w$ is differentiable at $z=1$ with $w\left( 1 \right) = 1$, then
$${\rm{Re}}\left[{{w_z}(1)+{w_{\bar z}}(1)}\right]\geq\frac{2}{\pi}\,\frac{(p+1)+\frac{\pi }{4}(p-1)(|a_p|+|b_p|)}{1+\frac{\pi}{4}(|a_p|+|b_p|)},$$
where $a_p=\frac{D^ph(0)}{p!}$ and $b_p=\frac{D^pg(0)}{p!}$.
\end{thm}

For $p=1$, it follows from (\ref{April-22}) that $|a_1|+|b_1|\leq\frac{4}{\pi}$. By using Theorem \ref{thm-1.1}, we then have
$${\rm{Re}}\left[{{w_z}(1)+{w_{\bar z}}(1)}\right]\geq\frac{2}{\pi}\,\frac{2}{1+\frac{\pi}{4}(|a_1|+|b_1|)}\geq\frac{2}{\pi}.$$

\begin{thm}\label{thm-1.2}
Let $w=h+\bar{g}$ be a sense-preserving harmonic mapping of $\ID$ such that $\mu(a, w)=p\geq1$ and $w(\ID)\subset\ID$, where $a\in \ID$.
If $w$ is differentiable at $z=\alpha$ with $w(\alpha)=\beta$, where $\alpha, \beta\in\mathbb{T}$, then
$${\rm{Re}}\bigg(\bar{\beta}[w_z(\alpha)\alpha+w_{\bar{z}}(\alpha)\bar{\alpha}]\bigg)\geq\frac{2}{\pi}\frac{(p+1)+\frac{\pi }{4}(p-1)\Lambda_w^{(p)}(a)(1-|a|^2)^p}{1+\frac{\pi}{4}\Lambda_w^{(p)}(a)(1-|a|^2)^p}\frac{1-|a|^2}{|1-\bar{a}\alpha|^2},$$
where $\Lambda_w^{(p)}(a)=\left|\frac{D^ph(a)}{p!}\right|+\left|\frac{D^pg(a)}{p!}\right|$.

In particular, when $\alpha=\beta=1$ and $a=0$, then Theorem $\ref{thm-1.2}$ coincides with Theorem $\ref{thm-1.1}$.
\end{thm}

The rest of this paper is organized as follows: in Section \ref{sec-2} we shall introduce some known results and prove two lemmas which will be used in the proof of our main results; in Section \ref{sec-3} we should prove Lemma \ref{bhz-lem-1}, Theorem \ref{thm-1.1} and Theorem \ref{thm-1.2}.
\section{Auxiliary results}\label{sec-2}
The following lemmas will be used in proving our main results.

\begin{lem}\label{HLZ-lemma1}\cite[Theorem 2]{Roman}
If $m(t)$ and $q(t)$ are functions for which all the necessary derivatives are defined, then
$${D^n}m(q(t))=\sum\limits_{{k_1}+\cdots+ n{k_n}=n} {\frac{n!}{{k_1}! \cdots {k_n}!}}(D^{{k_1}+ \cdots +{k_n}}m)(q(t))\left(\frac{D(q(t))}{1!}\right)^{k_1} \cdots\left(\frac{{D^n}(q(t))}{n!}\right)^{k_n},$$
where $k_1, \cdots, k_n$ are non-negative integer numbers.
\end{lem}

\begin{lem}\label{HLZ-lemma2}
Let $S = \left\{ {w\in\mathbb{C}:\;| {{\mathop{\rm Re}\nolimits} (w)}| < 1} \right\}$ be a strip domain, and $f: \ID\rightarrow S$ be an analytic function such that $\mu(0, f)= p \ge 1$. Assume that $\delta(z)=\tan \left( {\frac{\pi}{4}f(z)} \right)$. Then $\delta(z)$ is analytic in $\ID$ with  $\delta(\ID)\subset\ID$ and $\mu(0, \delta)=p\ge 1$.
\end{lem}
\bpf
We first prove that $\delta(z)$ is analytic in $\ID$ and $\delta(\ID)\subset\ID$.

To show this, assume that $f(z)=u+iv$ and let
$$\zeta=e^{\frac{\pi f}{2}i}=e^{-\frac{\pi v}{2}}e^{\frac{\pi u}{2}i}.$$
Since $f$ is an analytic function of $\ID$ into $S$, we see that $\zeta$ is an analytic function of $\ID$ into
$\mathbb{H}_+=\{\zeta\in\mathbb{C}:\; {\rm Re} \zeta>0\}$. This implies that $\delta(z)$ is analytic in $\ID$, since
$$\delta(z)=(-i)\frac{\zeta(z)-1}{\zeta(z)+1}.$$
The M\"obius transformation $\frac{\zeta-1}{\zeta+1}$ maps $\mathbb{H}_+$ into $\ID$, and thus, $\delta(\ID)\subset\ID$.

Secondly, we show that $\mu(0, \delta)=p$.

Obviously, $\delta(0)=0$. Let $\varphi(z)=\frac{\pi}{4}f(z)$. It follows from Lemma \ref{HLZ-lemma1} that

\begin{align*}
D^n\delta&=\sum\limits_{k_1 + \cdots  + n{k_n} = n}\frac{n!}{k_1!\cdots k_n!} \left(D^{k_1+\cdots +k_n}\tan\right)(\varphi)\left(\frac{D \varphi}{1!}\right)^{k_1}\cdots\left(\frac{D^n \varphi}{n!}\right)^{k_n}\\
&=\sum\limits_{\substack{k_1 + \cdots  + n{k_n} = n \\ \ \; k_n=0}}\frac{n!}{k_1!\cdots k_n!} \left(D^{k_1+\cdots +k_n}\tan\right)(\varphi)\left(\frac{D \varphi}{1!}\right)^{k_1}\cdots\left(\frac{D^n \varphi}{n!}\right)^{k_n}\\
&+(D\tan)(\varphi)D^n\varphi.
\end{align*}
The condition $\mu(0, f)=p$ ensures that $D^k \varphi(0)=0$, for $k=1, 2, \cdots , p-1$ and $D^p \varphi(0)\neq0$.
Therefore, $D^n \delta(0)=0$, for $n=1, 2, \cdots, p-1$. For $n=p$, we have
\be\label{eq2.1}D^p \delta(0)=(D\tan)(\varphi(0))D^p\varphi(0)=D^p\varphi(0)\neq0,\ee
which shows that $\mu(0, \delta)=p$.
\epf

Given $a\in\ID$, let $\eta(z)=\varphi _a(z)= \frac{a - z}{1 -{\bar{a}z}} $ be an automorphism of $\ID$, which interchanges $a$ and $z$.
Then we have the following lemma.
\begin{lem}\label{HLZ-lemma3}
Let $w=h+\bar{g}$ be a sense-preserving harmonic mapping of $\ID$ such that $\mu(a,w)=p \ge1$ and $w(\ID) \subset \ID$, where $a\in\ID$.
Assume that $W= w\circ{\eta}$. Then $W$ is a sense-preserving harmonic self-mapping of $\ID$ and $\mu(0,W)= p.$
\end{lem}
\bpf
Obviously,
$$W(0)=w(\eta(0))=w(a)=0.$$
Elementary calculations show that,
$$|W_z(z)|-|W_{\bar{z}}(z)|=(|w_\eta(\eta)|-|w_{\bar{\eta}}(\eta)|)|\varphi'_a(z)|.$$
Since $w$ is sense-preserving in $\ID$, we see that $|W_z|-|W_{\bar{z}}|>0$, and thus, $W$ is also sense-preserving in $\ID$.
Using Lemma \ref{HLZ-lemma1}, we obtain
\begin{align}\label{eq2.2}\nonumber
{D^n}W
&=\sum\limits_{{k_1}+\cdots+n{k_n}=n}\frac{n!}{{k_1}!\cdots{k_n}!}(D^{{k_1}+\cdots+{k_n}}w)(\eta){\left(\frac{D\eta}{1!}\right)}^{k_1} \cdots {\left(\frac{D^n\eta}{n!}\right)}^{k_n}\\
&=\sum\limits_{\substack{{k_1}+\cdots+n{k_n}=n \\ \;{k_1}\neq n}}\frac{n!}{{k_1}!\cdots{k_n}!}(D^{{k_1}+\cdots+{k_n}}w)(\eta){\left(\frac{D\eta}{1!}\right)}^{k_1} \cdots {\left(\frac{D^n\eta}{n!}\right)}^{k_n}\\ \nonumber
&+({D^n}w)(\eta){(D\eta)}^n.
\end{align}
The condition $\mu(a,w)=p$  ensures that
$$D^kw(a)=D^kw(\eta(0))=0, \ \ \  \mbox{where}\ \ \ k=1,2,\cdots,p-1.$$
Note that in (\ref{eq2.2}), if $k_1\neq n$, then $k_1+\cdots+k_n<n$, and thus $(D^{k_1+\cdots+k_n}w)(\eta(0))=0$.
Then
\be\label{March-31-1}D^nW(0)=0,\ \ \ \mbox{where}\ \ \ n=1,2,\cdots,p-1.\ee

For $n=p$, we have
\be\label{April-1-1}D^pW(0)=(D^pw)(\eta(0)){(\eta'(0))}^p.\ee
Since $\eta'(0)=|a|^2-1$ and
$$(D^pw)(\eta(0))=D^pw(a)\ne 0,$$
we see that
\be\label{March-31-2}D^pW(0)\neq0.\ee
Hence, $\mu(0,W)=p$ easily follows from (\ref{March-31-1}) and (\ref{March-31-2}).
\epf

\section{Main results}\label{sec-3}
\subsection{Proof of Lemma \ref{bhz-lem-1}.}
Assume that $w=u+iv$ is a sense-preserving harmonic self-mapping of $\ID$ with $\mu(0, w)=p\geq1$.
For any $\theta\in[0, 2\pi]$, let $f$ be an analytic function of $\ID$, where
$$\mbox{Re}f=u\cos\theta+v\sin\theta$$
is harmonic in $\ID$.
Then $f(\ID)\subset S=\{z\in\mathbb{C}:\;|{\rm Re}z|<1\}$ and $f(0)=0$.
If we write $f=\xi+i\vartheta$ and $w=h+\bar{g}$, then for $z=x+iy\in\ID$,
$$\xi(z)={\rm Re}(w(z)e^{-i\theta}),$$
and
\begin{align*}
f'(z)&=\xi_x(z)-i\xi_y(z) \\
&=h'(z)e^{-i\theta}+g'(z)e^{i\theta}.
\end{align*}
Therefore
\be\label{eq3.1}
{D^p}f={D^p}he^{-i\theta}+{D^p}ge^{i\theta},
\ee
which shows that $\mu(0,f)=p$, since $\mu(0,w)=p$.
Let
$$\delta=\tan\left(\frac{\pi}{4}f\right).$$
Then by Lemma \ref{HLZ-lemma2}, we see that $\delta$ is an analytic function of $\ID$ into $\ID$
with $\mu(0,\delta)= p$. Applying Lemma \Ref{HLZ-lemma4}, we have
$$|\delta(z)|\leq {|z|}^p\frac{|z|+\frac{1}{p!}|{D^p}\delta(0)|}{1+\frac{1}{p!}|{D^p}\delta(0)||z|}={|z|}^p\frac{|z|+\frac{\pi}{4p!}|{D^p}f(0)|}{1+\frac{\pi}{4p!}|{D^p}f(0)||z|},$$
where the last equality holds since it follows from (\ref{eq2.1}) that $D^p\delta(0)=\frac{\pi}{4}D^pf(0)$.

On the other hand, let
$$d(z)=\frac{e^{i\frac{\pi}{2}f(z)}-1}{e^{i\frac{\pi}{2}f(z)}+1}.$$
Then $d(z)=i\delta(z)$.
Using the following elementary inequality
$$\tan\frac{1}{2}|{\rm Re}\varsigma|\leq\left|\frac{e^{i\varsigma}-1}{e^{i\varsigma}+1}\right|,\ \ \ \mbox{for all}\ \ \ |{\rm Re}\varsigma|\leq\frac{\pi}{2},$$
we see that $$\tan\left(\frac{1}{2}\left|{\rm Re}\frac{\pi}{2}f\right|\right)\leq |d|=|\delta|.$$

Thus
\be\label{March-31-4}|{\rm Re}f(z)|\leq\frac{4}{\pi}\arctan|\delta|\leq\frac{4}{\pi}\arctan
\left[{|z|}^p\frac{|z|+\frac{\pi}{4p!}|{D^p}f(0)|}{1+\frac{\pi}{4p!}|{D^p}f(0)||z|}\right].\ee

Using (\ref{eq3.1}), we have
$$\left|\frac{{D^p}f(0)}{p!}\right|=\left|\frac{{D^p}h(0)e^{-i\theta}}{p!}+\frac{{D^p}g(0)e^{i\theta}}{p!}\right|\leq\left|\frac{{D^p}h(0)}{p!}\right|+\left|\frac{{D^p}g(0)}{p!}\right|=|a_p|+|b_p|.$$
Elementary calculations show that for $0 \leq r < 1$, the function $\varphi \left( x \right) = \frac{{r + \frac{\pi }{4}x}}{{1 + \frac{\pi }{4}xr}}$ is an increasing function of $x$. These together with (\ref{March-31-4}) show that
\be\label{March-31-3}|u(z)\cos\theta+v(z)\sin\theta|\leq\frac{4}{\pi}\arctan
\left[{|z|}^p\frac{|z|+\frac{\pi}{4}(|a_p|+|b_p|)}{1+\frac{\pi}{4}(|a_p|+|b_p|)|z|}\right].\ee

The desired inequality (\ref{bhz-lem1-1}) is now easy to follow, since
$$|w(z)|=\max\limits_{\theta\in[0, 2\pi]}|\xi|=\max\limits_{\theta\in[0, 2\pi]}|u(z)\cos\theta+v(z)\sin\theta|.$$
This completes the proof of Lemma \ref{bhz-lem-1}.
\qed

\subsection{Proof of Theorem \ref{thm-1.1}.}
For any $z\in\ID$, since $\mu(0, w)=p$, we see from Lemma \ref{bhz-lem-1} that
\be\label{eq3.2}
|w(z)|\leq\frac{4}{\pi}\arctan
\left[{|z|}^p\frac{|z|+\frac{\pi}{4}(|a_p|+|b_p|)}{1+\frac{\pi}{4}(|a_p|+|b_p|)|z|}\right]:=M(|z|).
\ee
Since $w$ is differential at $z=1$, we know that
$$w(z)=1+{w_z}(1)(z-1)+w_{\bar{z}}(1)(\bar{z}-1)+\circ(|z-1|).$$
This together with (\ref{eq3.2}) show that
$$|1+{w_z}(1)(z-1)+{w_{\bar{z}}}(1)(\bar{z}-1)+\circ(|z-1|)|^2 \le M^2(|z|).$$
Therefore,
$$2{\rm Re}[{w_z}(1)(1-z)+w_{\bar{z}}(1)(1-\bar{z})] \ge 1-M^2(|z|)+ \circ(|z-1|).$$
Take $z=r\in(0,1)$ and letting $r \to {1^-}$, it follows from $M(1)=1$ that
\begin{align*}
 2{\rm{Re}}[{w_z}(1)+{w_{\bar{z}}}(1)] &\ge \mathop {\lim }\limits_{r \to {1^-}} \frac{1-M^2(r)}{1-r} \\
 &= \frac{4}{\pi}\frac{(p+1)+\frac{\pi}{4}(p-1){(|a_p|+|b_p|)}}{{1 + \frac{\pi}{4}{(|a_p|+|b_p|)}}}.
\end{align*}
Then
\be
{\rm{Re}}[{w_z}(1)+{w_{\bar{z}}}(1)] \ge \frac{2}{\pi}\frac{(p+1)+\frac{\pi}{4}(p-1){(|a_p|+|b_p|)}}{1+\frac{\pi}{4}(|a_p|+|b_p|)},
\ee
hence the proof of the theorem is complete.
\qed

\subsection{Proof of Theorem \ref{thm-1.2}.}
For $\alpha\in\mathbb{T}$, let $\gamma= {\eta}(\alpha)$, where $\eta(z)=\varphi _a(z)= \frac{a - z}{1 -{\bar{a}z}}$. It is easy to see that $\gamma\in\mathbb{T}$ and
$${\eta}(\gamma)=\alpha. $$
Elementary calculations show that
$$\eta'(0)=|a|^2- 1.$$

For $\beta \in \mathbb{T}$, let $W(\zeta)= \bar{\beta} w \circ {\eta}(\zeta \gamma)=H+\bar{G}$, where $\zeta \in \ID$.
Then
$${W_\zeta}(\zeta)=\bar{\beta}{w_z}(\eta(\zeta \gamma)){\eta'}(\zeta \gamma)\gamma $$
and
$$W_{\bar{\zeta}}(\zeta)=\bar{\beta}{w_{\bar z }}({\eta}(\zeta \gamma))\overline{{\eta'}(\zeta \gamma)} \bar{\gamma}.$$
Using the following equation
$$\eta'(\gamma)=\frac{-(1-\bar{a}\alpha)^2}{1-|a|^2}$$
we have
\be\label{April-1-2}
{\rm Re}({W_\zeta}(1)+{W_{\bar{\zeta}}}(1))={\rm Re}\left(\bar{\beta}\left[{w_z}(\alpha)\alpha \frac{{|1 - \bar{a}\alpha|}^2}{1-{|a|}^2}+w_{\bar{z}}(\alpha)\bar{\alpha}\frac{{|1-\bar{a}\alpha |}^2}{1-{|a|}^2}\right]\right).
\ee
Since $w$ is a sense-preserving harmonic self-mapping of $\ID$ with $\mu(a, w)=p$, it follows from Lemma \ref{HLZ-lemma3} that
$W(\zeta)$ is also sense-preserving in $\ID$ with $W(\ID)\subset\ID$ and $\mu(0,W)=p$.
Furthermore, we have
$$W(0)=\bar {\beta}w(\eta(0))=\bar{\beta}w(a)=0$$
and
$$W(1)=\bar{\beta}w(\eta(\gamma))=\bar{\beta}w(\alpha)=|\beta|^2=1.$$
Using Theorem \ref{thm-1.1}, we obtain the following inequality
\be\label{April-1-3}
{\rm Re}({W_\zeta}(1)+{W_{\bar{\zeta}}}(1))\ge \frac{2}{\pi}\frac{(p+1)+\frac{\pi }{4}(p-1)\left(\frac{1}{p!}|{D^p}H(0)|+ \frac{1}{p!}|{D^p}G(0)|\right)}{1 + \frac{\pi }{4}\left(\frac{1}{p!}|{D^P}H(0)| +\frac{1}{p!}|{D^P}G(0)|\right)}.
\ee
According to (\ref{April-1-1}) and note that $\bar{D}^pW(0)=(\bar{D}^pw)(\eta(0))\left(\bar{\eta}'(0)\right)^p$, we have
\be\label{April-1-4}
\frac{1}{p!}|{D^p}H(0)|+ \frac{1}{p!}|{D^p}G(0)|=\Lambda_w^{(p)}(a)(1-|a|^2)^p,\ee
where $\Lambda_w^{(p)}(a)=\left|\frac{D^pw(a)}{p!}\right|+\left|\frac{\bar{D}^pw(a)}{p!}\right|$.
It follows from (\ref{April-1-2}), (\ref{April-1-3}) and (\ref{April-1-4}) that
$${\rm Re}\bigg(\bar{\beta}[{w_z}(\alpha)\alpha+{w_{\bar{z}}}(\alpha)\bar{\alpha}]\bigg) \ge \frac{2}{\pi}\frac{(p+1)+\frac{\pi }{4}(p-1){\Lambda_w^{(p)}}(a)(1-|a|^2)^p}{1+ \frac{\pi }{4}{\Lambda_w^{(p)}}(a)(1-|a|^2)^p}\frac{1-|a|^2}{|1-\overline a \alpha|^2}.$$
If $a = 0$, then
$${\rm Re}\bigg(\bar{\beta}[{w_z}(\alpha)\alpha+{w_{\bar {z}}}(\alpha)\bar{\alpha}]\bigg)\geq\frac{2}{\pi}\frac{(p+1)+\frac{\pi }{4}(p-1)(|a_p|+|b_p|)}{1+\frac{\pi}{4}(|a_p|+|b_p|)}.$$
This completes the proof of the theorem.
\qed

\vspace*{5mm}
\noindent {\bf Acknowledgments}. We would like to thank the anonymous referees for their helpful comments to improve this paper.

\vspace*{5mm}
\noindent {\bf Funding}. The research of the authors were supported by
NNSF of China Grant Nos. 11471128, 11501220, 11971124, 11971182, NNSF of Fujian Province Grant Nos. 2016J01020, 2019J0101,  Subsidized Project for Postgraduates'Innovative Fund in Scientific Research of Huaqiao University and the Promotion Program for Young and Middle-aged Teacher in Science and Technology Research of Huaqiao University (ZQN-PY402).

\end{document}